\documentclass[11pt]{amsart}
\usepackage{amsmath}
\usepackage{amsthm}
\usepackage{latexsym}
\usepackage{amssymb}
\usepackage{xspace}
\usepackage{amscd}

\newcommand{\Gal}{{\rm Gal}}

\newcommand{\Aut}{\mbox{\rm Aut}}

\newcommand{\Ker}{\mbox{\rm Ker}}

\newcommand{\inv}{^{-1}}

\newcommand{\Z}{{\mathbb Z}}
\newcommand{\Q}{{\mathbb Q}}

\newcommand{\matriz}[1]{\begin{array} #1 \end{array}}

\newcommand{\GEN}[1]{\langle #1 \rangle}

\newcommand{\wh}[1]{\widehat{#1}}

\newcommand{\G}{\overline{G}}

\title{The Schur group of an abelian number field}

\author[A. Herman]{Allen Herman}
\address{Department of Mathematics, University of Regina, Regina, Canada}
\email{aherman@math.uregina.ca}
\author[G. Olteanu]{Gabriela Olteanu}
\address{Department of Mathematics and Computer Science, North University of Baia Mare,
Victoriei 76, 430122 Baia Mare, Romania.
Current address. Departamento de Matem\'{a}ticas,
Universidad de Murcia, Murcia, Espa\~{n}a}
\email{olteanu@math.ubbcluj.ro, golteanu@um.es}
\author[\'{A}. del R\'{\i}o]{\'{A}ngel del R\'{\i}o}
\address{Departamento de Matem\'{a}ticas, Universidad de Murcia, Murcia, Espa\~{n}a}
\email{adelrio@um.es}

\thanks{{Research supported by the National Science and Engineering Research
Council of Canada, D.G.I. of Spain and Fundaci\'{o}n
S\'{e}neca of Murcia.}}

\newtheorem{theorem}{Theorem}
\newtheorem{lemma}[theorem]{Lemma}
\newtheorem{proposition}[theorem]{Proposition}

\newtheorem{corollary}[theorem]{Corollary}

\newtheorem{notation}[theorem]{Notation}

\theoremstyle{remark}
\theoremstyle{remark}

\begin{document}

\begin{abstract}
We characterize the maximum $r$-local index of a Schur algebra over an abelian number field $K$ in terms of global information determined by the field $K$, for $r$ an arbitrary rational prime. This completes and unifies previous results of Janusz in \cite{Jan} and Pendergrass in \cite{Pen1}.
\end{abstract}

\maketitle

\section{Introduction and Preliminaries}

Let $K$ be a field. A {\em Schur algebra} over $K$ is a central simple $K$-algebra which is generated over $K$ by a finite group of units. The {\em Schur group} of $K$ is the subgroup $S(K)$ of the Brauer group of $K$ formed by classes containing a Schur algebra.  By the Brauer-Witt Theorem (see e.g. \cite{Yam}), each class in $S(K)$ can be represented by a cyclotomic algebra, i.e. a crossed product of the form $(L/K, \alpha)$ in which $L/K$ is a cyclotomic extension and the factor set $\alpha$ takes values in the group of roots of unity $W(L)$ of $L$.

In the case when $K$ is an abelian number field; i.e. $K$ is contained in a finite cyclotomic extension of $\Q$, Benard-Schacher theory \cite{BS} gives a partial characterization of the elements of $S(K)$.  According to this theory, if $n$ is the Schur index of a Schur algebra over $K$, then $W(K)$ contains an element of order $n$. This is known as the Benard-Schacher Theorem. Furthermore, if $\frac{t}{n}$ (in lowest terms) is the local invariant of $A$ at a prime $\mathcal{R}$ of $K$ that lies over a rational prime $r$, then each of the fractions $\frac{c}{n}$ with $1 \le c \le n$ and $c$ coprime to $n$ will occur equally often among the local invariants corresponding to the primes of $K$ lying above $r$.  In particular, these local invariants all have the same denominator $n$ for all the primes of $K$ lying above $r$, which we call the {\it $r$-local index} $m_r(A)$ of $A$. Only finitely many of the $m_r(A)$ are greater than $1$, and the Schur index of $A$ is the least common multiple of the $m_r(A)$ as $r$ runs over all rational primes.

The goal of this article is to characterize the maximum $r$-local index of a Schur algebra over an abelian number field $K$ in terms of global information determined by $K$. The existence of this maximum is a consequence of the Benard-Schacher Theorem. Since $S(K)$ is a torsion abelian group, it is enough to compute the maximum of the $r$-local indices of Schur algebras over $K$ with index a power of $p$ for every prime $p$ dividing the order of $W(K)$.  We will refer to this number as $p^{\beta_p(r)}$.  In \cite{Jan}, Janusz gave a formula for $p^{\beta_p(r)}$ when either $p$ is odd or $K$ contains a primitive $4$-th root of unity. The remaining cases were considered by Pendergrass in \cite{Pen1}. However, some of the calculations involving factor sets in \cite{Pen1} are not correct, and as a consequence the formulas for $2^{\beta_2(r)}$ for odd primes $r$ that appear there are inaccurate. This article was motivated in part to find a correct formula for $p^{\beta_p(r)}$ in this remaining case, and also because of the need to apply the formula in an upcoming work of the authors in \cite{HOR}.  Since the local index at $\infty$ will be $2$ when $K$ is real and will be $1$ otherwise, the only  remaining case is that of $r=2$.  In this case, $p$ must be equal to $2$ and we must have $\zeta_4 \not\in K$. The characterization of fields $K$ for which $S(K_2)$ is of order $2$ is given in \cite[Corollary 3.3]{Pen1}.

The main result of the paper (Theorem \ref{localindex}) characterizes $p^{\beta_p(r)}$ in terms of the position of $K$ relative to an overlying cyclotomic extension $F$ that is determined by $K$ and $p$.  The formulas for $p^{\beta_p(r)}$ are stated in terms of elements of certain Galois groups in this setting.  The main difference between our approach and that of Janusz and Pendergrass is that the field $F$ that we use is slightly larger, which allows us to present some of the somewhat artificial-looking calculations in \cite{Jan} in a more conceptual fashion.  Another highlight of our approach is the treatment of calculations involving factor sets.  In Section~\ref{SecFSC} we generalize a result from \cite{AS} which describes the factor sets for a given action of an abelian group $G$ on another abelian group $W$ in terms of some data.  In particular, we give necessary and sufficient conditions that the data must satisfy in order to be induced by a factor set.  Because of the applications we have in mind, extra attention is paid to the case when $W$ is a cyclic $p$-group.

\section{Factor set calculations}\label{SecFSC}

In this section $W$ and $G$ are two abelian groups and $\Upsilon:G\rightarrow \Aut(W)$ is a group homomorphism. A group epimorphism $\pi:\G\rightarrow G$ with kernel $W$ is said to induce $\Upsilon$ if, given $u_g\in \G$ such that $\pi(u_g)=g$, one has $u_g w u_g\inv=\Upsilon(g)(w)$ for each $w\in W$.  If $g \mapsto u_g$ is a {\it crossed section} of $\pi$ (i.e. $\pi(u_g)=g$ for each $g \in G$) then the map $\alpha: G \times G \rightarrow W$ defined by $u_g u_h= \alpha_{g,h} u_{gh}$ is a {\it factor set} (or {\it $2$-cocycle}) $\alpha \in Z^2(G,W)$.  We always assume that the crossed sections are normalized, i.e. $u_1=1$ and hence $\alpha_{g,1}=\alpha_{1,g}=1$.  Since a different choice of crossed section for $\pi$ would be a map $g \mapsto w_g u_g$ where $w: G \rightarrow W$, $\pi$ determines a unique cohomology class in $H^2(G,W)$, namely the one represented by $\alpha$.

Given a list $g_1,\dots,g_n$ of generating elements of $G$, a group epimorphism $\pi:\G \rightarrow G$ inducing $\Upsilon$,
and a crossed section $g\mapsto u_g$ of $\pi$, we associate the elements $\beta_{ij}$ and $\gamma_i$ of $W$, for $i,j\le
n$, by the equalities:
    \begin{equation}\label{BetaGamma}
    \matriz{{rcl} u_{g_j}u_{g_i} &=& \beta_{ij} u_{g_i} u_{g_j}, \mbox{ and } \\
    u_{g_i}^{q_i} & = &\gamma_i u_{g_1}^{t_1^{(i)}}\cdots u_{g_{i-1}}^{t_{i-1}^{(i)}},}
    \end{equation}
where the integers $q_i$ and $t_j^{(i)}$ for $1\le i \le n$ and $0\le j < i$ are determined by
    \begin{equation}\label{qt}
    q_i = \text{ order of } g_i \text{ modulo } \langle g_1,\dots,g_{i-1} \rangle, \quad
    g_i^{q_i}=g_1^{t_1^{(i)}}\cdots g_{i-1}^{t_{i-1}^{(i)}}, \quad \text{and} \quad
    0\le t_j^{(i)} < q_j.
    \end{equation}
If $\alpha$ is the factor set associated to $\pi$ and the crossed section $g\mapsto u_g$, then we say that $\alpha$
induces the data $(\beta_{ij},\gamma_i)$. The following proposition gives necessary and sufficient conditions for a
list $(\beta_{ij},\gamma_i)$ of elements of $W$ to be induced by a factor set.

The order of an element $g$ of a group is denoted by $|g|$.

\begin{proposition}\label{Ext}
Let $W$ and $G=\GEN{g_1,\dots,g_n}$ be abelian groups and let $\Upsilon:G\rightarrow \Aut(W)$ be an action of $G$ on $W$. For
every $1 \le i,j \le n$,
let $q_i$ and $t_j^{(i)}$ be the integers determined by
(\ref{qt}).  For every $w\in W$ and $1 \le i \le n$, let
$$\Upsilon_i=\Upsilon(g_i), \quad N_i^t(w) = w \Upsilon_i(w) \Upsilon_i^2(w) \cdots \Upsilon_i^{t-1}(w), \quad \text{and} \quad N_i = N_i^{q_i}.$$

For every $1\le i,j \le n$, let $\beta_{ij}$ and $\gamma_i$ be elements of $W$. Then the following conditions are
equivalent:
\begin{enumerate}
\item There is a factor set $\alpha\in Z^2(G,W)$ inducing the data $(\beta_{ij},\gamma_i)$.
\item The following equalities hold for every $1\le i,j,k \le n$:
\begin{itemize}
\item[(C1)] $\beta_{ii}=\beta_{ij}\beta_{ji}=1$.
\item[(C2)] $\beta_{ij} \beta_{jk} \beta_{ki}  = \Upsilon_k(\beta_{ij}) \Upsilon_i(\beta_{jk}) \Upsilon_j(\beta_{ki})$.
\item[(C3)] $N_i(\beta_{ij})\gamma_i =
    \Upsilon_j(\gamma_i) N_1^{t_1^{(i)}}(\beta_{1j}) \Upsilon_1^{t_1^{(i)}}(N_2^{t_2^{(i)}}(\beta_{2j})) \cdots
    \Upsilon_1^{t_1^{(i)}}\Upsilon_2^{t_2^{(i)}} \dots \Upsilon_{i-2}^{t_{i-2}^{(i)}}
    (N_{i-1}^{t_{i-1}^{(i)}}(\beta_{(i-1)j}))$.
\end{itemize}
\end{enumerate}
\end{proposition}

\begin{proof}
(1) implies (2). Assume that there is a factor set $\alpha\in Z^2(G,W)$ inducing the data $(\beta_{ij},\gamma_i)$. Then
there is a surjective homomorphism $\pi:\G\rightarrow G$ and a crossed section $g \mapsto u_g$ of $\pi$ such that
the $\beta_{ij}$ and $\gamma_i$ satisfy (\ref{BetaGamma}). Condition (C1) is clear. Conjugating by $u_{g_k}$ in
$u_{g_j}u_{g_i} = \beta_{ij} u_{g_i} u_{g_j}$ yields
    $$\matriz{{c}
    \beta_{jk} \Upsilon_j(\beta_{ik}) \beta_{ij} u_{g_i} u_{g_j} = \beta_{jk} \Upsilon_j(\beta_{ik}) u_{g_j} u_{g_i} = \beta_{jk} u_{g_j} \beta_{ik} u_{g_i} =
    u_{g_k} u_{g_j} u_{g_i} u_{g_k}\inv = \\ u_{g_k} \beta_{ij}  u_{g_i} u_{g_j} u_{g_k}\inv =
    \Upsilon_k(\beta_{ij}) \beta_{ik}  u_{g_i} \beta_{jk} u_{g_j}
    = \Upsilon_k(\beta_{ij}) \beta_{ik} \Upsilon_i(\beta_{jk}) u_{g_i} u_{g_j}.}$$
Therefore, we have $\beta_{jk} \Upsilon_j(\beta_{ik}) \beta_{ij} = \Upsilon_k(\beta_{ij}) \beta_{ik}
\Upsilon_i(\beta_{jk})$ and so (C2) follows from (C1).

To prove (C3), we use the obvious relation $(wu_{g_i})^t = N_i^t(w)u_{g_i}^t$.  Conjugating by $u_{g_j}$ in
    $u_{g_i}^{q_i}=\gamma_iu_{g_1}^{t_1^{(i)}}\cdots u_{g_{i-1}}^{t_{i-1}^{(i)}}$ results in
    $$\matriz{{c}
    N_i(\beta_{ij})\gamma_i u_{g_1}^{t_1^{(i)}} \cdots u_{g_{i-1}}^{t_{i-1}^{(i)}} = N_i^{q_i}(\beta_{ij})u_{g_i}^{q_i}  = (\beta_{ij}u_{g_i})^{q_i} =
    u_{g_j} u_{g_i}^{q_i} u_{g_j}\inv = u_{g_j} \gamma_i u_{g_1}^{t_1^{(i)}} \cdots u_{g_{i-1}}^{t_{i-1}^{(i)}} u_{g_j}\inv = \\
    \Upsilon_j(\gamma_i) (\beta_{1j} u_{g_1})^{t_1^{(i)}} \cdots (\beta_{(i-1)j} u_{g_{i-1}})^{t_{i-1}^{(i)}} =
    \Upsilon_j(\gamma_i) N_1^{t_1^{(i)}}(\beta_{1j}) u_{g_1}^{t_1^{(i)}} \cdots
        N_{i-1}^{t_{i-1}^{(i)}}(\beta_{(i-1)j}) u_{g_{i-1}}^{t_{i-1}^{(i)}} = \\
    \Upsilon_j(\gamma_i) N_1^{t_1^{(i)}}(\beta_{1j}) \Upsilon_1^{t_1^{(i)}}(N_2^{t_2^{(i)}}(\beta_{2j})) \cdots
        \Upsilon_1^{t_1^{(i)}}\Upsilon_2^{t_2^{(i)}} \dots \Upsilon_{i-2}^{t_{i-2}^{(i)}} (N_{i-1}^{t_{i-1}^{(i)}}(\beta_{(i-1)j}))
        u_{g_1}^{t_1^{(i)}} \cdots u_{g_{i-1}}^{t_{i-1}^{(i)}}.
    }$$
Cancelling on both sides produces (C3). This finishes the proof of (1) implies (2).

Before proving (2) implies (1), we show that if $\pi:\G\rightarrow G$ is a group homomorphism with kernel $W$ inducing
$\Upsilon$, $g \mapsto u_g$ is a crossed section of $\pi$ and $\beta_{ij}$ and $\gamma_i$ are given by
(\ref{BetaGamma}), then $\G$ is isomorphic to the group $\wh{G}$ given by the following presentation: the set of generators of $\wh{G}$
is $\{\wh{w},\wh{g}_i:w\in W, i=1,\dots,n\}$, and the relations are
    \begin{equation}\label{GenRel}
    \wh{w_1w_2}=\wh{w_1}\wh{w_2},\quad \Upsilon_i(w)=\wh{g}_i \wh{w} \wh{g}_i\inv, \quad \wh{g}_j\wh{g}_i = \wh{\beta}_{ij} \wh{g}_i \wh{g}_j \quad \text{ and }\quad
    \wh{g}_i^{q_i}=\wh{\gamma_i} \wh{g}_1^{t_1^{(i)}}\cdots \wh{g}_{i-1}^{t_{i-1}^{(i)}},
    \end{equation}
for each $1\le i,j \le n$ and $w,w_1,w_2\in W$.
Since the relations obtained by replacing $\wh{w}$ by $w$ and $\wh{g}_i$ by $u_{g_i}$ in equation (\ref{GenRel})
for each $x\in W$ and each $1\le i \le n$, hold in $\G$, there is a surjective group homomorphism $\phi:\wh{G}
\rightarrow \G$, which associates $\wh{w}$ with $w$, for every $w\in W$, and $\wh{g}_i$ with $u_{g_i}$, for every
$i=1,\dots,n$. Moreover, $\phi$ restricts to an isomorphism $\wh{W} \rightarrow W$ and
$|\wh{g}_i\GEN{\wh{W},\wh{g}_1,\dots,\wh{g}_{i-1}}|=q_i$. Hence $[\wh{G}:\wh{W}|=q_1\cdots q_n=[\overline{G}:W]$ and so
$|\wh{G}|=|\overline{G}|$. We conclude that $\phi$ is an isomorphism.

(2) implies (1). Assume that the $\beta_{ij}$'s and $\gamma_i$'s satisfy conditions (C1), (C2) and (C3). We will recursively construct groups $\G_0,\G_1,\dots,\G_n$. Start with $\G_0=W$. Assume that $\G_{k-1}=\GEN{W,u_{g_1},\dots,u_{g_{k-1}}}$
has been constructed with $u_{g_1}, \dots, u_{g_{k-1}}$ satisfying the last three relations of (\ref{GenRel}), for $1\le i,j < k$, and that
these relations, together with the relations in $W$, form a complete list of relations for $\G_{k-1}$. To define $\G_k$
we first construct a semidirect product $H_k=\G_{k-1}\rtimes_{c_k} \GEN{x_k}$, where $c_k$ acts on $\G_{k-1}$ by
    $$c_k(w) = \Upsilon_k(w), \quad (w\in W), \quad \quad c_k(u_{g_i}) = \beta_{ik} u_{g_i}.$$
In order to check that this defines an automorphism of $\G_{k-1}$ we need to check that $c_k$ respects the defining
relations of $\G_{k-1}$. This follows from the commutativity of $G$ and conditions (C1), (C2) and (C3) by
straightforward calculations which we leave to the reader.

Notice that the defining relations of $H_k$ are the defining relations of $\G_{k-1}$ and the relations
$x_kw=\Upsilon_k(w)x_k$ and $x_ku_{g_i} = \beta_{ik} u_{g_i} x_k$. Using (C3) one deduces
    $u_{g_i} x_k^{q_k} u_{g_i}\inv = u_{g_i} \gamma_k u_{g_1}^{t_1^{(k)}} \cdots u_{g_{k-1}}^{t_{k-1}^{(k)}}
u_{g_i}\inv$, for each $i\le k-1$. This shows that $y_k=x_k^{-q_k} \gamma_k u_{g_1}^{t_1^{(k)}} \cdots
u_{g_{k-1}}^{t_{k-1}^{(k)}}$ belongs to the center of $H_k$. Let $\G_k=H_k/\GEN{y_k}$ and $u_{g_k} = x_k \GEN{ y_k }$. Now it is easy to see that
the defining relations of $G_k$ are the relations of $W$ and the last three relations in (\ref{GenRel}), for $0\le i,j
\le k$.

It is clear now that the assignment $w\mapsto 1$ and $u_{g_i} \mapsto g_i$ for each $i=1,\dots,n$ defines a group
homomorphism $\pi:\G=\G_n\rightarrow G$ with kernel $W$ and inducing $\Upsilon$. If $\alpha$ is the factor set
associated to $\pi$ and the crossed section $g\mapsto u_g$, then $(\beta_{ij},\gamma_i)$ is the list of data
induced by $\alpha$.
\end{proof}

Note that the group generated by the values of the factor set $\alpha$ coincides with the group generated by the
data $(\beta_{ij},\gamma_i)$. This observation will be used in the next section.

In the case $G=\GEN{g_1}\times \dots \times \GEN{g_n}$ we obtain the following corollary that one should compare with
Theorem 1.3 of \cite{AS}.

\begin{corollary}
If $G=\GEN{g_1}\times \dots \times \GEN{g_n}$ then a list $D=(\beta_{ij},\gamma_i)_{1\le i,j \le n}$ of elements of $W$
is the list of data associated to a factor set in $Z^2(G,W)$ if and only if the elements of $D$ satisfy {\rm (C1)},
{\rm (C2)} and $N_i(\beta_{ij})\gamma_i = \Upsilon_j(\gamma_i)$, for every $1\le i,j \le n$.
\end{corollary}

In the remainder of this section we assume that $W=\GEN{\zeta}$ is a cyclic $p$-group, for $p$ a prime integer.
Let $p^a$ and $p^{a+b}$ denote the orders of $W^G=\{x\in W : \Upsilon(g)(x) =x \text{ for each } g\in G\}$ and $W$
respectively. We assume that $0<a,b$. We also set
    $$C =  \Ker(\Upsilon) \quad \text{and} \quad D = \{g\in G : \Upsilon(g)(\zeta)=\zeta \text{ or }
    \Upsilon(g)(\zeta)=\zeta\inv\}.$$

Note that $D$ is subgroup of $G$ containing $C$, $G/D$ is cyclic, and $[D:C]\le 2$. Furthermore, the assumption $a>0$ implies that if $C \ne D$ then $p^a=2$.

\begin{lemma}\label{rho}
There exists a $\rho\in D$ and a subgroup $B$ of $C$ such that $D=\GEN{\rho}\times B$ and
$C=\GEN{\rho^2}\times B$.
\end{lemma}

\begin{proof}
The lemma is obvious if $C=D$ (just take $\rho=1$). So assume that $C \ne D$ and temporarily
take $\rho$ to be any element of $D\setminus C$. Since $[D:C]=2$, one may assume without loss of generality that $|\rho|$ is a power of $2$. Write $C=C_2\times C_{2'}$, where $C_2$ and $C_{2'}$ denote the $2$-primary and $2'$-primary parts of $C$, and
choose a decomposition $C_2=\GEN{c_1}\times \dots \times \GEN{c_n}$ of $C_2$. By reordering the $c_i$'s if needed, one
may assume that $\rho^2=c_1^{a_1}\dots c_k^{a_k} c_{k+1}^{2c_{k+1}} \dots c_n^{2a_n}$ with $a_1,\dots,a_k$ odd. Then
replacing $\rho$ by $\rho c_{k+1}^{-a_{k+1}}\dots c_n^{-a_n}$ one may assume that $\rho^2=c_1^{a_1}\dots c_k^{a_k}$,
with $a_1,\dots,a_k$ odd. Let $H=\GEN{\rho,c_1,\dots,c_k}$. Then $|\rho|/2=|\rho^2|=\exp(H\cap C)$, the exponent of
$H\cap C$, and so $\rho$ is an element of maximal order in $H$. This implies that $H=\GEN{\rho}\times H_1$ for some
$H_1\le H$. Moreover, if $h\in H_1\setminus C$ then $1\ne \rho^{|\rho|/2}=h^{|\rho|/2} \in \GEN{\rho}\cap H_1$, a
contradiction. This shows that $H_1\subseteq C$. Thus $C_2=(H\cap C_2) \times \GEN{c_{k+1}}\times \dots \times
\GEN{c_n} = \GEN{\rho^2} \times H_1 \times \GEN{c_{k+1}}\times \dots \times \GEN{c_n}$. Then $\rho$ and $B=H_1 \times
\GEN{c_{k+1}}\times \dots \times \GEN{c_n} \times C_{2'}$ satisfy the required conditions.
\end{proof}

By Lemma~\ref{rho}, there is a decomposition $D=B\times\GEN{\rho}$ with $C=B\times \GEN{\rho^2}$, which will be fixed
for the remainder of this section. Moreover, if $C=D$ then we assume $\rho=1$. Since $G/D$ is cyclic, $G/C=\GEN{\rho C}
\times \GEN{\sigma C}$ for some $\sigma\in G$.  It is easy to see that $\sigma$ can be selected so that if $D=G$
then $\sigma=1$, and $\sigma(\zeta)=\zeta^c$ for some integer $c$ satisfying
\begin{equation}\label{cee}
v_p(c^{q_{\sigma}}-1) = a + b, \mbox{ and }
v_p(c-1) = \left\{\matriz{{cl} a & \mbox{ if } G/C \mbox{ is cyclic and } G \ne D, \\
                             a+b & \mbox{ if } G/C \mbox{ is cyclic and $G = D$, and } \\
                         d \ge 2 & \mbox{ for some integer $d$, if } G/C \mbox{ is not cyclic, }}\right.
\end{equation}
where $q_{\sigma}=|\sigma C|$ and the map $v_p:\Q\rightarrow \Z$ is the classical $p$-adic valuation. In
particular, if $G/C$ is non-cyclic (equivalently $C\ne D\ne G$) then $p^a=2$, $b\ge 2$, $\rho(\zeta)=\zeta\inv$ and
$\sigma(\zeta^{2^{b-1}})=\zeta^{2^{b-1}}$.

For every positive integer $t$ we set
$$V(t) = 1+c+c^2+\dots+c^{t-1} = \frac{c^t-1}{c-1}.$$

Now we choose a
decomposition $B=\GEN{c_1}\times \dots \times \GEN{c_n}$ and adapt the notation of Proposition~\ref{Ext} for a group epimorphism
$f:\G\rightarrow G$ with kernel $W$ inducing $\Upsilon$ and elements
$u_{c_1},\dots,u_{c_n},u_{\sigma},u_{\rho}\in \G$ with $f(u_{c_i})=c_i$, $f(u_{\rho})=\rho$ and $f(u_{\sigma})=\sigma$,
by setting
    $$\beta_{ij}=[u_{c_j},u_{c_i}], \quad
    \beta_{i\rho}=\beta_{\rho i}\inv = [u_{\rho},u_{c_i}], \quad \beta_{i\sigma}=\beta_{\sigma i}\inv = [u_{\sigma},u_{c_i}], \mbox{ and }
     \beta_{\sigma\rho}=\beta_{\rho\sigma}\inv = [\beta_{\rho},\beta_{\sigma}].$$
We also set
    \begin{equation}\label{qtc}
    \matriz{{c}
    q_i=|c_i|, \quad q_{\rho}=|\rho|, \quad
    \text{and} \quad \sigma^{q_{\sigma}} = c_1^{t_1}\dots c_n^{t_n} \rho^{2t_{\rho}}, \\
    \text{ where } 0\le t_i < q_i \text{ and } 0\le t_{\rho} < |\rho^2|.}
    \end{equation}

With a slightly different notation than in Proposition~\ref{Ext}, we have, for each $1\le i \le n$, $t_j^{(i)}=0$ for
each $0\le j < i$, $t_i^{(\rho)}=0$ , $t_{i}^{(\sigma)}=t_i$, and $t_{\rho}^{(\sigma)}=2t_{\rho}$. Furthermore,
$q_{\rho}=1$ if $C=D$ and $q_{\rho}$ is even if $C\ne D$. Continuing with the adaptation of the notation of
Proposition~\ref{Ext} we set
    $$\gamma_i = u_{c_i}^{q_i}, \quad \gamma_{\rho}=u_{\rho}^{q_{\rho}}, \mbox{ and }
    \gamma_{\sigma}=u_{\sigma}^{q_{\sigma}}u_{c_1}^{-t_1}\dots u_{c_n}^{-t_n}u_{\rho}^{2t_{\rho}}.$$

We refer to the list $\{\beta_{ij},\beta_{i\sigma},\beta_{i\rho}, \beta_{\sigma\rho}, \gamma_i,\gamma_{\rho},
\gamma_{\sigma}: 0\le i<j\le n\}$, which we abbreviate as $(\beta, \gamma)$, as the data associated to the group epimorphism $f:\G\rightarrow G$ and choice of crossed section
$u_{c_1},\dots,u_{c_n},u_{\sigma},u_{\rho}$, or as the data induced by the corresponding factor set in $Z^2(G,W)$.

Furthermore, for every $w\in W$, $1 \le i \le n$ and $t\ge 0$ one has
    $$N_i^t(w)=w^t, \quad N_{\sigma}^t(w)=w^{V(t)}\quad \mbox{and} \quad
    N_{\rho}^t(w) = \left\{\matriz{{ll}
    w^t, & \mbox{if } \rho=1; \\
    1, & \mbox{if } \rho\ne 1 \text{ and } t \text{ is even};\\
    w, & \mbox{if } \rho\ne 1 \text{ and } t \text{ is odd}.}\right.
    $$
In particular, for every $w\in W$ one has
$$N_i(w)=w^{q_i}, \quad N_{\sigma}(w)=w^{V(q_{\sigma})}, \quad \text{and} \quad N_{\rho}(w)=1.$$

Rewriting Proposition~\ref{Ext} for this case we obtain the following.

\begin{corollary}\label{ExtCyc}
Let $W$ be a finite cyclic $p$-group and let $G$ be an abelian group acting on $W$ with
$G=\GEN{c_1,\dots,c_n,\sigma,\rho}$, $B=\GEN{c_1}\times \dots \times \GEN{c_n}$, $D=B\times \GEN{\rho}$ and $C=B\times
\GEN{\rho^2}$ as above. Let $q_i,q_{\rho},q_{\sigma}$ and the $t_i$'s be given by (\ref{qtc}). Let
$\beta_{\sigma\rho},\gamma_{\rho},\gamma_{\sigma}\in W$ and for every $1\le i,j \le n$ let
$\beta_{ij},\beta_{i\sigma},\beta_{i\rho}$ and $\gamma_i$ be elements of $W$. Then the following conditions are
equivalent:
\begin{enumerate}
\item The given collection $(\beta,\gamma) = \{\beta_{ij}, \gamma_i, \beta_{i\sigma}, \gamma_{\sigma}, \gamma_{\rho}, \beta_{\sigma\rho} \}$ is the list of data
induced by some factor set in $Z^2(G,W)$.
\item The following equalities hold for every $1\le i,j \le n$:
\begin{itemize}
\item[(C1)] $\beta_{ii}=\beta_{ij}\beta_{ji}=1$.
\item[(C2)]
\begin{enumerate}
\item $\beta_{ij}\in W^G$.
\item If $\rho\ne 1$ then $\beta_{i\sigma}^2=\beta_{i\rho}^{1-c}$.
\end{enumerate}
\item[(C3)]
    \begin{enumerate}
    \item $\beta_{ij}^{q_i}=1$.
    \item $\beta_{i\sigma}^{q_i}=\gamma_i^{c-1}$.
    \item $\beta_{i\sigma}^{-V(q_{\sigma})}=\beta_{1i}^{t_1}\dots \beta_{ni}^{t_n}$.
    \item $\gamma_{\sigma}^{c-1}\beta_{1\sigma}^{t_1}\dots \beta_{n\sigma}^{t_n} = 1$.
    \item If $\rho=1$ then $\beta_{i\rho}=\beta_{\sigma\rho}=\gamma_{\rho}=1$.
    \item If $\rho\ne 1$ then $\beta_{i\rho}^{q_i}\gamma_i^2=1$,
        $\beta_{\sigma\rho}^{V(q_{\sigma})}\gamma_{\sigma}^2=\beta_{1\rho}^{t_1}\dots \beta_{n\rho}^{t_n}$ and
        $\gamma_{\rho}\in W^G$.
    \end{enumerate}
\end{itemize}
\end{enumerate}
\end{corollary}

\begin{proof}
By completing the data with $\beta_{\sigma i}=\beta_{i\sigma}\inv$, $\beta_{\rho i}=\beta_{i\rho}\inv$ and
$\beta_{\sigma\sigma}=\beta_{\rho\rho}=1$ we have that (C1) is a rewriting of condition (C1) from
Proposition~\ref{Ext}.

(C2) is the rewriting of condition (C2) from Proposition~\ref{Ext} because this condition vanishes when $1 \le i,j,k \le n$
and when two of the elements $i,j,k$ are equal. Furthermore, permuting $i,j,k$ in (C2) yields equivalent conditions. So
we only have to consider three cases: substituting $i=i$, $j=j$, and $k=\sigma$;  $i=i$, $j=j$, and  $k=\rho$; and $i=i$, $j=\rho$, and $k=\sigma$.  In the first two cases one obtains $\sigma(\beta_{ij})=\rho(\beta_{ij})=\beta_{ij}$, or equivalently
$\beta_{ij}\in W^G$. For $\rho=1$ the last case vanishes, and for $\rho\ne 1$ (C2) yields
$\beta_{i\sigma}^2=\beta_{i\rho}^{1-c}$.

Rewriting (C3) from Proposition~\ref{Ext} we obtain: (a) for $i=i$, $j=j$; (b) for $i=i$ and $j=\sigma$; (c) for $i=\sigma$ and $j = i$; and (d) for $i=\sigma$ and $j=\sigma$.

We consider separately the cases $\rho=1$ and $\rho\ne 1$ for the remaining cases for rewriting (C3).  Assume first that $\rho=1$. When $i$ is replaced by $\rho$ and $j$ replaced by $i$ (respectively, by $\sigma$) we obtain
$\beta_{i\rho}=1$ (respectively $\beta_{\sigma\rho}=1$). On the other hand the requirement of only using normalized
crossed sections implies $\gamma_{\rho}=1$ in this case.  When $j=\rho$ the conditions obtained are trivial.

Now assume that $\rho\ne 1$.  For $i=i$ and $j=\rho$ one obtains $\beta_{i\rho}^{q_i}\gamma_i^2=1$. For $i=\rho$ and
$j=i$ one obtains a trivial condition because $N_{\rho}(x)=1$. For $i=\sigma$ and $j=\rho$, we
obtain $\beta_{\sigma\rho}^{V(q_{\sigma})}\gamma_{\sigma}^2=\beta_{1\rho}^{t_1}\dots \beta_{n\rho}^{t_n}$. For $i=\rho$ and
$j=\sigma$ one has $\sigma(\gamma_{\rho})=\gamma_{\rho}$, and for $i=\rho$ and $j=\rho$ one obtains
$\rho(\gamma_{\rho})=\gamma_{\rho}$.  The last two equalities are equivalent to $\gamma_{\rho}\in W^G$.
\end{proof}

\begin{corollary}\label{Sign}
With the notation of Corollary~\ref{ExtCyc}, assume that $G/C$ is non-cyclic
and $q_k$ and $t_k$
are even for some $k\le n$. Let $(\beta, \gamma)$ be the list of data induced by
a factor set
in $Z^2(G,W)$. Then the list obtained by replacing $\beta_{k\sigma}$ by $-\beta_{k\sigma}$ and keeping the remaining
data fixed is also induced by a factor set
in $Z^2(G,W)$.
\end{corollary}

\begin{proof}
It is enough to show that $\beta_{k\sigma}$ appears in all the conditions of Corollary~\ref{ExtCyc} with an even
exponent. Indeed, it only appears in (C2.b) with exponent 2; in (C3.b) with exponent $q_k$; in (C3.c) with exponent
$-V(q_{\sigma})$; and in (C3.d) and (C3.f) with exponent $t_k$. By the assumption it only remains to show that $V(q_{\sigma})$ is even.  Indeed, $v_2(V(q_{\sigma}))=v_2(c^{q_{\sigma}}-1)-v_2(c-1)=1+b-v_2(c-1)\ge 1$ because $c \not\equiv 1 \mod 2^{1+b}$.
\end{proof}

The data $(\beta, \gamma)$ induced by
a factor set are not cohomologically invariant because they depend on the selection of $\pi$ and of the $u_{c_i}$'s, $u_{\sigma}$ and $u_{\rho}$. However, at least the $\beta_{ij}$ are cohomologically invariant. For every $\alpha\in H^2(G,W)$ we associate a matrix $\beta_{\alpha}=(\beta_{ij})_{1\le i,j \le n}$ of
elements of $W^G$ as follows: First select a group epimorphism $\pi:\G \rightarrow G$ realizing $\alpha$
and $u_{c_1},\dots,u_{c_n}\in \G$ such that $\pi(u_{c_i})=c_i$, and then
set $\beta_{ij}=[u_{c_j},u_{c_i}]$. The definition of $\beta_{\alpha}$ does not depend on the choice of $\pi$ and the
$u_{c_i}$'s because  if $w_1,w_2\in W$ and $u_1,u_2\in \G$ then $[w_1u_1,w_2u_2]=[u_1,u_2]$.

\begin{proposition}\label{SP}
Let $\beta=(\beta_{ij})_{1\le i,j \le n}$ be a matrix of elements of $W^G$ and for every $1\le i,j\le n$ let
$a_{ii}=0$ and $a_{ij}=\min(a,v_p(q_i),v_p(q_j))$, if $i\ne j$.

Then there is an $\alpha\in H^2(G,W)$ such that $\beta=\beta_{\alpha}$ if and only if the following conditions hold for every $1\le i,j \le n$:
\begin{equation}\label{beta}
\beta_{ij}\beta_{ji}=\beta_{ij}^{p^{a_{ij}}}=1.
\end{equation}
\end{proposition}

\begin{proof}
Assume first that $\beta=\beta_{\alpha}$ for some $\alpha\in Z^2(G,W)$. Then (\ref{beta}) is a consequence of
conditions (C1), (C2.a) and (C3.a) of Corollary~\ref{ExtCyc}.

Conversely, assume that $\beta$ satisfies (\ref{beta}). The idea of the proof is that one can enlarge $\beta$ to
a list of data $(\beta, \gamma)$
that satisfies conditions (C1)--(C3) of Corollary~\ref{ExtCyc}. Hence the desired conclusion follows
from the corollary.

Condition (C1) follows automatically from (\ref{beta}). If
$i,j\le n$ then $\beta_{ij}\in W^G$ follows from the fact that $a\ge a_{ij}$ and so (\ref{beta}) implies that
$\beta_{ij}^{p^a}=1$. Hence (C2.a) holds. Also (C3.a) holds automatically from (\ref{beta}) because $p^{a_{ij}}$
divides $q_i$. Hence, we have to select the $\beta_{i\sigma}$'s, $\beta_{i\rho}$'s, $\gamma_i$'s, $\beta_{\sigma \rho}$, $\gamma_{\sigma}$, and $\gamma_{\rho}$  for (C2.b) and (C3.b)--(C3.f) to hold.

Assume first that $D=G$.
In this case we just take
$\beta_{i{\sigma}}=\beta_{i{\rho}}=\beta_{\sigma\rho}=\gamma_i=\gamma_{\sigma}=\gamma_{\rho}=1$ for every $i$. Then
(C2.b), (C3.b), (C3.d) and (C3.f) hold trivially by our selection. Moreover, in this case $\sigma=1$ and so $t_i=0$ for
each $i=1,\dots,n$, hence (C3.c) also holds.

In the remainder of the proof we assume that $D\ne G$.  First we show how one can assign values to
$\beta_{{\sigma}i}$ and $\gamma_i$, for $i\le n$ for (C3.b)--(C3.d) to hold.  Let $d=v_p(c-1)$ and
$e=v_p(V(q_{\sigma}))=a+b-d$. (see (\ref{cee})).
Note that $d=a$ if $C=D$ and $a=1\le 2 \le d \le b$ if $C\ne D$ (because we
are assuming that $D\ne G$). Let $X_1,X_2,Y_1$ and $Y_2$ be integers such that $c-1=p^dX_1$, $V(q_{\sigma})=p^eX_2$,
and $X_1Y_1\equiv X_2Y_2 \equiv 1 \mod p^{a+b}$. By (\ref{beta}), $\beta_{ij}^{p^{a_{ij}}}=1$ and so $\beta_{ij}\in
W^{p^{a+b-a_{ij}}}$. Therefore there are integers $b_{ij}$, for $1\le i,j\le n$ such that $b_{ii}=b_{ij}+b_{ij}=0$
and $\beta_{ij}=\zeta^{b_{ij}p^{a+b-a_{ij}}}$. For every $i\le n$ set
    $$x_i=Y_2\sum_{j=1}^n t_j b_{ji} p^{a-a_{ji}}, \quad
    \beta_{\sigma i} = \zeta^{x_ip^{d-a}} \quad
    y_i=Y_1Y_2\sum_{j=1}^n t_j b_{ji} \frac{q_i}{p^{a_{ij}}}, \quad \text{and} \quad
    \gamma_i = \zeta^{y_i}.$$
Then $V(q_{\sigma})p^{d-a}x_i = p^eX_2Y_2 \sum_{j=1}^n t_j b_{ji} p^{d-a_{ji}} \equiv \sum_{j=1}^n t_j b_{ji} p^{a+b-a_{ji}}
\mod p^{a+b}$ and therefore
    $$\beta_{{\sigma}i}^{V(q_{\sigma})}= \zeta^{\sum_{j=1}^{n} t_j b_{ji} p^{a+b-a_{ji}}}= \prod_{i=1}^{n} \beta_{ji}^{t_j},$$
that is (C3.c) holds. Moreover $q_i p^{d-a} x_i = p^d Y_2 \sum_{j=1}^n t_j b_{ji} \frac{q_i}{p^{a_{ij}}} \equiv p^d X_1
y_i = (c-1)y_i$ and therefore $\beta_{i\sigma}^{q_i} = \gamma_i^{c-1}$, that is (C3.b) holds.

We now compute
    \begin{equation}\label{Sum}
    \sum_{i=1}^n t_i x_i = Y_2\sum_{1\le i,j \le n} t_i t_j b_{ij} p^{a-a_{ij}} =
    Y_2\sum_{i=1}^{n+1} t_i^2 b_{ii} p^{a-a_{ii}} + Y_2\sum_{1\le i<j \le n} t_i t_j (b_{ij}+b_{ji}) p^{a-a_{ij}}
    = 0.
    \end{equation}
Then setting $\gamma_{\sigma}=1$, one has
    $$\gamma_{\sigma}^{c-1} \prod_{i=1}^{n} \beta_{i{\sigma}}^{t_i} = \prod_{i=1}^{n} \zeta^{-t_ix_ip^{d-a}} =
    \zeta^{-p^{d-a}\sum_{i=1}^{n} t_ix_i} = 1 $$
and (C3.d) holds. This finishes the assignments of $\beta_{i\sigma}$ and $\gamma_i$ for $i\le n$ and of
$\gamma_{\sigma}$.

If $C=D$ then a quick end is obtained assigning $\beta_{i{\rho}}=\beta_{\sigma\rho}=\gamma_{\rho}=1$.

So it only remains to assign values to $\beta_{i{\rho}}, \beta_{\sigma\rho}$ and $\gamma_{\rho}$ under the
assumption that $C\ne D$.  Set $\beta_{i\rho}=\zeta^{-Y_1x_i}$. In this case $p^a=2$ and therefore $2p^{d-a}x_i = p^d
x_i \equiv (c-1) Y_1 x_i$ and $q_iY_1x_i = 2y_i$. Thus $\beta_{i\sigma}^2\beta_{i\rho}^{c-1} = \zeta^{2p^{d-a}x_i}
\zeta^{(1-c)Y_1x_i}=1$, hence (C2.b) holds, and $\beta_{i\rho}^{q_i}\gamma_i^2 = \zeta^{-q_iY_1x_i+2y_i}=1$, hence the
first relation of (C3.f) follows.

Finally, using (\ref{Sum}) one has
    $$\beta_{1\rho}^{t_1}\dots \beta_{n\rho}^{t_n} = (\beta_{1\sigma}^{t_1}\dots \beta_{n\sigma}^{t_n})^{-Y_1} = 1 =
    \gamma_{\sigma}^2$$
and the last two relations of (C3.f) hold when $\beta_{\sigma\rho}=\gamma_{\rho}=1$.
\end{proof}

Let $\beta=(\beta_{ij})$ be an $n\times n$ matrix of elements of $W^G$ satisfying (\ref{beta}). Then the map
$\Psi:B\times B \rightarrow W^G$ given by
    $$\Psi((c_1^{x_1} \dots c_n^{x_n},c_1^{y_1} \dots c_n^{y_n})) = \prod_{1\le i,j \le n}
    \beta_{ij}^{x_iy_j}$$
is a {\it skew pairing} of $B$ over $W^G$ in the sense of \cite{Jan}; that is, it satisfies the following conditions
for every $x,y,z \in B$:
    $$(\Psi 1) \quad \Psi(x,x)=\Psi(x,y)\Psi(y,x)=1,
    \quad\quad (\Psi 2) \quad \Psi(x,yz)=\Psi(x,y)\Psi(x,z).$$
Conversely, every skew pairing of $B$ over $W^G$ is given by a matrix $\beta=(\beta_{ij}=\Psi(c_i,c_j))_{1\le i,j
\le n}$ satisfying (\ref{beta}). In particular, every class in $H^2(G,W)$ induces a skew pairing $\Psi=\Psi_{\alpha}$ of $B$ over $W^G$ given by $\Psi(x,y)=\alpha_{x,y}\alpha_{y,x}^{-1}$, for all $x,y \in B$, for any cocycle $\alpha$ representing the given cohomology class.

In terms of skew pairings, Proposition~\ref{SP} takes the following form.

\begin{corollary}\label{FacSets}
If $\Psi$ is a skew pairing of $B$ over $W^G$ then there is an $\alpha\in H^2(G,W)$ such that $\Psi=\Psi_{\alpha}$.
\end{corollary}

Corollary~\ref{FacSets} was obtained in \cite[Proposition 2.5]{Jan} for $p^a\ne 2$. The remaining cases were considered
in \cite[Corollary 1.3]{Pen1}, where it is stated that for every skew pairing $\Psi$ of $C$ over $W^G$ there is a
factor set $\alpha\in Z^2(G,W)$ such that $\Psi(x,y)=\alpha_{x,y}\alpha_{y,x}\inv$, for all $x,y \in C$. However,
this is false if $\rho^2\ne 1$ and $B$ has nontrivial elements of order $2$.  Indeed, if $\Psi$ is the skew pairing of
$B$ over $W^G$ given by the factor set $\alpha$ then $\Psi(x,\rho^2)=1$ for each $x \in C$. To see this we introduce a
new set of generators of $G$, namely $G=\GEN{c_1,\dots,c_n,c_{n+1},\rho,\sigma}$ with $c_{n+1}=\rho^2$. Then condition
(C3) of Proposition~\ref{Ext}, for $i=\rho$ and $j=i$ reads $\beta_{(n+1)i}=1$ which is equivalent to
$\Psi(c_i,\rho^2)=1$ for all $1 \le i \le n$. Using this it is easy to give a counterexample to \cite[Corollary
1.3]{Pen1}.

Before finishing this section we mention two lemmas that will be needed in next section.  The first is elementary and so the proof has been omitted.

\begin{lemma}\label{Maximums}
Let $S$ be the set of skew pairings of $B$ with values in $W^G$.
If $B=B' \times B''$ and $b_1, b_2 \in B'$ and $b_3 \in B''$ then $$\max\{\Psi(b_1 \cdot b_3,b_2): \Psi \in S\}= \max\{\Psi(b_1,b_2):\Psi \in S\} \cdot \max\{\Psi(b_3,b_2):\Psi \in S\}.$$
\end{lemma}

\begin{lemma}\label{SPEx}
Let $\widehat{B} = B \times \GEN{g}$ be an abelian group and let $h\in B$. If $k=\gcd\{p^a,|g|\}$ and $t=|hB^k|$ then $t$ is the maximum possible value of $\Psi(h,g)$ as $\Psi$ runs over all skew pairings of $\widehat{B}$ over $\GEN{ \zeta_{p^a} }$.
\end{lemma}

\begin{proof}
Since $k$ divides $p^a$, the hypothesis $t=|hB^k|$ implies that there is a group homomorphism $\chi:B\rightarrow \GEN{
\zeta_{p^a}}$ such that $\chi(B^k)=1$ and $\chi(h)$ has order $t$. Let $\Psi:\widehat{B}\times \widehat{B} \rightarrow \GEN{
\zeta_{p^a} }$ be given by $\Psi(xg^i,yg^j)=\chi(x^jy^{-i})=\chi(x)^i\chi(y)^{-j}$, for $x,y\in B$. If
$g^i=g^{i'}$, then $i\equiv i' \mod |g|$ and hence $i\equiv i' \mod k$.  Therefore, $x^i B^k = x^{i'} B^k$, which implies that
$\chi(x)^i=\chi(x)^{i'}$. This shows that $\Psi$ is well defined. Now it is easy to see that $\Psi$ is a skew pairing
and $\Psi(h,g)=\chi(h)$ has order $t$.

Conversely, if $\Psi$ is any skew pairing of $\widehat{B}$ over $\GEN{ \zeta_{p^a} }$, then $\Psi(x,g)^{p^a} = 1$ and
$\Psi(x,g)^{|x|} = \Psi(1,g) = 1$ for all $x \in B$.  This implies that $\Psi(x^k,g)=\Psi(x,g)^k=1$ for all $x \in B$,
and so $\Psi(B^k,g)=1$. Therefore $\Psi(h,g)^t = \Psi(h^t,g) \in \Psi(B^k,g) = 1$, so the order of $\Psi(h,g)$ divides
$t$.
\end{proof}

\section{Local index computations}

In this section $K$ denotes an abelian number field, $p$ a prime, and $r$ an odd prime. Our goal is to find a global formula for $\beta(r) = \beta_p(r)$, the maximum nonnegative integer for which $p^{\beta(r)}$ is the $r$-local index of a Schur algebra over $K$.

We are going to abuse the notation and denote by $K_r$ the completion of $K$ at a (any) prime of $K$ dividing $r$. If
$E/K$ is a finite Galois extension, one may assume that the prime of $E$ dividing $r$, used to compute $E_r$, divides
the prime of $K$ over $r$, used to compute $K_r$. We use the classical notation:
    $$\matriz{{rcl}
    e(E/K,r) & = & e(E_r/K_r) =  \text{ramification index of } E_r/K_r. \\
    f(E/K,r) & = & f(E_r/K_r) =  \text{residue degree of } E_r/K_r. \\
    m_r(A) &=& \text{Index of } K_r\otimes_K A, \text{ for a Schur algebra } A \text{ over } K.}
    $$
By Benard-Schacher Theory and because $E/K$ is a finite Galois extension, $e(E/K,r)$, $f(E/K,r)$ and $m_r(A)$ do not depend on the selection of the prime of $K$ dividing $r$ (see \cite{Serre} and \cite{BS}).  By the Benard-Schacher Theorem and because $|S(K_r)|$ divides $r-1$ \cite{Yam}, if either $\zeta_p\not\in K$ or $r\not\equiv 1 \mod p$ then $\beta(r)=0$. So to avoid trivialities we assume that $\zeta_p\in K$ and $r \equiv 1 \mod p$.

Suppose $K \subseteq F=\Q(\zeta_n)$ for some positive integer $n$ and let $n=r^{v_r(n)}n'$. Then $\Gal(F/\Q)$
contains a {\it canonical Frobenius automorphism at $r$} which is defined by
$\psi_r(\zeta_{r^{v_r(n)}})=\zeta_{r^{v_r(n)}}$ and $\psi_r(\zeta_{n'})=\zeta_{n'}^r$.  We can then define the {\it
canonical Frobenius automorphism at $r$ in $\Gal(F/K)$} as $\phi_r = \psi_r^{f(F/K,r)}$.  On the other hand, the {\it
inertia subgroup at $r$ in $\Gal(F/K)$} is by definition the subgroup of $\Gal(F/K)$ that acts as
$\Gal(F_r/K_r(\zeta_{n'}))$ in the completion at $r$.

We use the following notation.

\begin{notation}\label{Not} First we define some positive integers:

    $m=$ minimum even positive integer with $K\subseteq \Q(\zeta_m)$,

    $a=$ minimum positive integer with $\zeta_{p^a}\in K$,

    $s = v_p(m)$ and

    $$b = \left\{ \matriz{{ll}
            s, & \mbox{ if $p$ is odd or $\zeta_4 \in K$}, \\
            s+v_p([K\cap \Q(\zeta_{p^s}):\Q]) + 2, &
                                \mbox{if } \Gal(K(\zeta_{p^{2a+s}})/K) \mbox{ is not cyclic, and } \\

            s+1, & \mbox{otherwise.} } \right.$$
We also define
    $$L=\Q(\zeta_m), \quad \zeta=\zeta_{p^{a+b}}, \quad W = \GEN{\zeta}, \quad F = L(\zeta),$$
    $$G=\Gal(F/K), \quad C=\Gal(F/K(\zeta)), \quad
    \text{and} \quad D = \Gal(F/K(\zeta+\zeta\inv)).$$

Since $\zeta_p\in K$, the automorphism $\Upsilon:G \rightarrow \Aut(W)$ induced by the Galois action satisfies the
conditions of Section~\ref{SecFSC} and the notation is consistent. As in that section we fix elements $\rho$ and
$\sigma$ in $G$ and a subgroup $B = \GEN{c_1} \times \dots \times \GEN{c_n}$ of $C$ such that $D=B\times
\GEN{\rho}$, $C=B\times \GEN{\rho^2}$ and $G/C=\GEN{\rho C}\times \GEN{\sigma C}$. Furthermore,
$\sigma(\zeta)=\zeta^c$ for some integer $c$ chosen according to (\ref{cee}). Notice that by the choice of $b$, $G\ne
B$.

We also fix an odd prime $r$ and set
    $$e=e(K(\zeta_r)/K,r), \quad f = f(K/\Q,r) \quad \text{and} \quad \nu(r)=\max\{0,a+v_p(e)-v_p(r^f-1)\}.$$
Let $\phi\in G$ be the canonical Frobenius automorphism at $r$ in $G$,
and write
    $$\phi = \rho^{j'} \sigma^j \eta,\quad \text{ with } \eta \in B,\quad 0
    \le j' < |\rho| \quad \text{ and } \quad 0 \le j < |\sigma C|.$$

Let $q$ be an odd prime not dividing $m$.
Let
$G_q=\Gal(F(\zeta_q)/K)$, $C_q=\Gal(F(\zeta_{q})/K(\zeta))$ and
let $c_0$ denote a generator of $\Gal(F(\zeta_q)/F)$. Finally we
fix

$\theta=\theta_q$, a generator of the inertia group of $r$ in
$G_q$ and

$\phi_q = c_0^{s_0}\phi =c_0^{s_0}\eta \rho^{j'} \sigma^j=\eta_q
\rho^{j'} \sigma^j$, the canonical Frobenius automorphism at $r$ in $G_q$.
\end{notation}

Observe that we are considering $G$ as a subgroup of $G_q$ by identifying $G$ with $\Gal(F(\zeta_q)/K(\zeta_q))$. Again
the Galois action induces a homomorphism $\Upsilon_q:G_q\rightarrow \Aut(W)$ and $W^{G_q}=\GEN{\zeta_{p^a}}$.
So this action satisfies the conditions of Section~\ref{SecFSC} and we adapt the notation by settting
    $$B_q=\GEN{c_0}\times B, \quad
    C_q = \Gal(F(\zeta_q)/K(\zeta)) = \Ker(\Upsilon_q) \quad \text{and} \quad
    D_q=\Gal(F(\zeta_q)/K(\zeta+\zeta\inv)).$$
Notice that $C_q = \GEN{c_0} \times C = B_q \times \GEN{\rho^2}$  and $D_q=D\times \GEN{c_0}$. Hence $G/C\simeq
G_q/C_q$.

If $\Psi$ is a skew pairing of $B$ over $W^G$ then $\Psi$ has a
unique extension to a skew pairing $\Psi$ of $C$ over $W^G$ which
satisfies $\Psi(B,\rho^2)=\Psi(\rho^2,B)=1$. So we are going to
apply skew pairings of $B$ to pairs of elements in $C$ under the
assumption that we are using this extension.

Since $p\ne r$, $\theta\in C_q$.  Moreover, if $r=q$ then $\theta$ is a generator of $\Gal(F(\zeta_r)/F)$ and otherwise
$\theta\in C$.  Notice also that if $G/C$ is
non-cyclic then $p^a=2$ and $K\cap \Q(\zeta_{2^s})=\Q(\zeta_{2^d}+\zeta_{2^d}\inv)$, where $d = v_p(c-1)$, and so
$b=s+d$.

It follows from results of Janusz \cite[Proposition 3.2]{Jan} and
Pendergrass \cite[Theorem 1]{Pen2} that $p^{\beta(r)}$ always
occurs as the $r$-local index of a cyclotomic algebra of the
form $(L(\zeta_q)/L,\alpha)$ where $q$ is either $4$ or a prime
not dividing $m$ and $\alpha$ takes values in $W(L(\zeta_q))_p$, with the
possibility of $q=4$ occurring only in the case when $p^s=2$.
By inflating the factor set $\alpha$ to $F(\zeta_q)$ (which will be equal to $F$ when $p^s=2$), we have that $p^{\beta(r)} = m_r(A)$, where
\begin{equation}\label{Algebra}
\matriz{{l} A = (F(\zeta_q)/K,\alpha) \text{ (we also write $\alpha$ for the inflation),} \\
q \text{ is an odd prime not dividing } m, \text{ and } \\
\alpha \text{ takes values in $\GEN{\zeta_{p^4}}$ if $p^s = 2$ and in $\GEN{\zeta_{p^{s}}}$ otherwise. }}
\end{equation}
So it suffices to find a formula for the maximum $r$-local index of a
Schur algebra over $K$ of this form.

Write $A=\bigoplus_{g\in G_q} F(\zeta_q)u_g$, with $u_g\inv x u_g
= g(x)$ and $u_gu_h=\alpha_{g,h}u_{gh}$, for each $x\in F(\zeta_q)$
and $g,h\in G_q$.  After a diagonal change of basis one may assume
that if $g=c_0^{s_0}c_1^{s_1}\dots
c_{n}^{s_n}\rho^{s_{\rho}}\sigma^{s_{\sigma}}$ with $0\le s_i <
q_i=|c_i|$, $0\le s_{\rho} < |\rho|$ and $0\le
s_{\sigma}<q_{\sigma}=|\sigma C|$ then $u_g=u_{c_0}^{s_0}u_{c_1}^{s_1}\dots
u_{c_{n}}^{s_n}u_{\rho}^{s_{\rho}}u_{\sigma}^{s_{\sigma}}$.

It is well known (see \cite{Yam} and \cite[Theorem 1]{Jan}) that
\begin{equation}\label{mxi}
m_r(A)=|\xi|, \quad \text{where} \quad
\xi=\xi_{\alpha}=\left(\frac{\alpha_{\theta,\phi_q}}{\alpha_{\phi_q,\theta}}\right)^{r^{v_r(e)}} u_{\theta}^{r^{v_r(e)}(r^f-1)}.
\end{equation}
This can be slightly simplified as follows. If $r|e$ then $\GEN{\theta}$ has an element $\theta^k$ of order $r$. Since
$\theta$ fixes every root of unity of order coprime with $r$, necessarily $r^2$ divides $m$ and the fixed field of
$\theta^k$ in $L$ is $\Q(\zeta_{m/r})$. Then $K\subseteq \Q(\zeta_{m/r})$, contradicting the minimality of $m$. Thus
$r\nmid e$ and so
    \begin{equation}\label{xi}
    \xi=\frac{\alpha_{\theta,\phi_q}}{\alpha_{\phi_q,\theta}} u_{\theta}^{r^f-1} =
    \frac{\alpha_{\theta,\phi_q}}{\alpha_{\phi_q,\theta}} \gamma_{\theta}^{\frac{r^f-1}{e}} =
    [u_{\theta},u_{\phi_q}] \gamma_{\theta}^{\frac{r^f-1}{e}}, \text{ where } \gamma_{\theta} = u_{\theta}^e.
    \end{equation}

With our choice of the $\{u_g: g \in G_q \}$, we have
    $$[u_{\theta},u_{\phi_q}] = [u_{\theta},u_{\eta_q}u_{\rho}^{j'} u_{\sigma}^j ] =
    \Psi(\theta,\eta_q)[u_{\theta},u_{\rho}^{j'} u_{\sigma}^j],$$
where $\Psi=\Psi_{\alpha}$ is the skew pairing associated to $\alpha$.  Therefore,
    $$\xi = \xi_0 \Psi(\theta,\eta_q) \quad \mbox{ with } \quad
    \xi_0 = \xi_{0,\alpha}=[u_{\theta},u_{\rho}^{j'}  u_{\sigma}^j] \gamma_{\theta}^{\frac{r^f-1}{e}}.$$

Let  $(\beta,\gamma)$ be the data associated to the factor set $\alpha$ (relative to the set of generators $c_1,\dots,c_n,\rho,\sigma$).

\begin{lemma}\label{xi0}
Let $A=(F(\zeta_q)/K,\alpha)$ be a cyclotomic algebra satisfying the conditions of (\ref{Algebra}) and use the above notation.  Let $\theta = c_0^{s_0} c_1^{s_1} \cdots c_n^{s_n} \rho^{2s_{n+1}}$, with $0\le s_i < q_i$ for $0 \le i\le n$, and $0\le s_{n+1}
\le |\rho^2|$.
\begin{enumerate}
\item If $G/C$ is cyclic then $\xi_0^{p^{\nu(r)}}=1$.
\item Assume that $G/C$ is non cyclic and let $\mu_i = \beta_{i\rho}^{\frac{1-c}{2}}\beta_{i\sigma}\inv$.
Then $\mu_i=\pm 1$ and $\xi_0^{p^{\nu(r)}} = \prod_{i=0}^n \mu_{i}^{2^{\nu(r)} (j+j') s_i}$.
\end{enumerate}
\end{lemma}

\begin{proof}
For the sake of regularity we write $c_{n+1}=\rho^2$.
Since $e=|\theta|$, we have that $q_i$ divides $es_i$ for each
$i$. Furthermore, $v_p(e)$ is the maximum of the
$v_p\left(\frac{q_i}{\gcd(q_i,s_i)}\right)$ for $i=1,\dots,n$.
Then
$$v_p(e)-v_p(r^f-1)=\max\left\{v_p\left(\frac{q_i}{\gcd(q_i,s_i)(r^f-1)}\right),i=1,\dots,n\right\}.$$
Hence
    \begin{equation}\label{nuMin}
    \matriz{{rcl}
    \nu(r)&=&\max\{0,v_p(e)+a-v_p(r^f-1)\} \\
    &=& \min\left\{x\ge 0:p^a \text{ divides } p^x \cdot \frac{s_i(r^f-1)}{q_i}, \mbox{ for each }
    i=1,\dots,n\right\}.}
    \end{equation}

Now we compute $\gamma_{\theta}$ in terms of the previous expression of $\theta$. Set $v=u_{c_{n+1}}^{s_{n+1}}$ and
$y=u_{c_0}^{s_0}u_{c_1}^{s_1}\cdots u_{c_n}^{s_n}$.  Then
    $$u_{\theta}=yv = \gamma vy, \quad \text{with} \quad
    \gamma= \Psi(c_{n+1}^{s_{n+1}},c_0^{s_0}c_1^{s_1}\dots,c_n^{s_n}).$$
Thus $\gamma^e =
\Psi(c_{n+1}^{es_{n+1}},c_0^{s_0}c_1^{s_1}\dots,c_n^{s_n}) = 1$.
Using that $[y,\gamma]=1$, one easily proves by induction on $m$
that
    $$(yv)^m = \gamma^{\binom{m}{2}} y^m v^m.$$
Hence
    $$(yv)^e = \gamma^{\binom{e}{2}} y^e v^e = \gamma^{\binom{e}{2}} y^e u_{c_{n+1}}^{es_{n+1}} =
    \gamma^{\binom{e}{2}} y^e \gamma_{\rho}^{\frac{es_{n+1}}{q_{n+1}}},$$
and $\gamma^{\binom{e}{2}}=\pm 1$. (If $p$ or $e$ is odd then necessarily $\gamma^{\binom{e}{2}}=1$.) Now an easy
induction argument shows
$$\gamma_{\theta} = \mu \gamma_0^{\frac{es_0}{q_0}} \gamma_1^{\frac{es_1}{q_1}} \cdots \gamma_n^{\frac{es_n}{q_n}}
\gamma_{\rho}^{\frac{es_{n+1}}{q_{n+1}}}, \quad \text{ for some } \mu = \pm 1.$$

Note that $\nu(r)+v_p(r^f-1)-v_p(e)\ge a \ge 1$, by (\ref{nuMin}). Then
$\mu^{p^{\nu(r)}\frac{r^f-1}{e}}=\gamma_{\rho}^{p^{\nu(r)}\frac{r^f-1}{e}}=1$, because both $\mu$ and $\gamma_{\rho}$
are $\pm 1$, and they are 1 if $p$ is odd (see (C3.e) and (C3.f)). Thus
    \begin{equation}\label{epstheta}
    \gamma_{\theta}^{p^{\nu(r)}\frac{r^f-1}{e}} =
    \prod_{i=0}^n \gamma_i^{p^{\nu(r)}\frac{(r^f-1)s_i}{q_i}}
    \end{equation}

(1). Assume that $G/C$ is cyclic.  We have that $\rho=1$ and $v_p(c-1)=a$.  Note that the $\beta$'s and $\gamma$'s are $p^b$-th roots of unity by (\ref{Algebra}).

Let $Y$ be an integer satisfying $Y\frac{c-1}{p^a} \equiv 1 \mod p^b$.  Since $\phi_q=\sigma^j \eta_q$ with
$\eta_q\in C_q$, we have $r^f\equiv c^j \mod p^{a+b}$ and so $Y\frac{r^f-1}{p^a}= Y\frac{c-1}{p^a} \frac{c^j-1}{c-1}
\equiv V(j)\mod p^b$.  Then $\beta_{i\sigma}^{Y\frac{r^f-1}{p^a}} = \beta_{i\sigma}^{V(j)}$.

Using that $p^a$ divides $p^{\nu(r)}\frac{s_i(r^f-1)}{q_i}$ (see (\ref{nuMin})) and $Y\frac{(c-1)}{p^a} \equiv 1
\mod p^b$ we obtain
    $$
    \gamma_{i}^{p^{\nu(r)}\frac{s_i(r^f-1)}{q_i}} =
    (\gamma_{i}^{c-1})^{Y\frac{p^{\nu(r)}s_i(r^f-1)}{p^aq_i}}.
    $$

 Combining this with (C3.b) we have
    \begin{equation}\label{Factor=1}
    \matriz{{rcl}
    [u_{c_i}^{s_i},u_{\sigma}^j]^{p^{\nu(r)}} \gamma_{i}^{p^{\nu(r)}\frac{s_i(r^f-1)}{q_i}} &=&
    [u_{c_i},u_{\sigma}]^{s_iV(j)p^{\nu(r)}} (\gamma_{i}^{c-1})^{Y\frac{p^{\nu(r)}s_i(r^f-1)}{p^aq_i}} \\
    &=& [u_{c_i},u_{\sigma}]^{s_iV(j)p^{\nu(r)}} \beta_{i\sigma}^{Y\frac{p^{\nu(r)}s_i(r^f-1)}{p^a}} \\
    &=& ([u_{c_i},u_{\sigma}]\beta_{i\sigma})^{p^{\nu(r)}s_iV(j)} = 1,}
    \end{equation}
because $\beta_{i\sigma} =[u_{\sigma},u_{c_i}] =
[u_{c_i},u_{\sigma}]\inv$. Using (\ref{epstheta}) and (\ref{Factor=1}) we have
    $$\xi_0^{p^{\nu(r)}} = [u_{\theta},u_{\sigma}^j]^{p^{\nu(r)}} \gamma_{\theta}^{p^{\nu(r)}\frac{r^f-1}{e}}
    = \prod_{i=0}^n [u_{c_i}^{s_i},u_{\sigma}^j]^{p^{\nu(r)}}
        \gamma_{i}^{p^{\nu(r)}\frac{s_i(r^f-1)}{q_i}} = 1$$
and the lemma is proved in this case.

(2).  Assume now that $G/C$ is non-cyclic. Then $p^a=2$ and if $d=v_2(c-1)$ then $d\ge 2$ and $b=s+d$. The data for
$\alpha$ lie in $\GEN{ \zeta_{2^{s+1}}}\subseteq \GEN{\zeta_{2^b}} \subseteq \GEN{\zeta_{2^{1+s+d}}} = W(F)_2$.  (C2.b) implies $\mu_i=\pm 1$ and using (C3.b) and (C3.f) one has $\gamma_{i}^{c+1}=\beta_{i\sigma}^{q_i}\beta_{i\rho}^{-q_i}$.  Let $X$ and $Y$ be integers satisfying $X\frac{c-1}{2^d} \equiv Y\frac{c+1}{2}\equiv 1\mod 2^{1+s+d}$ and set $Z=Y\frac{r^f-1}{2}$.

Recall that $2^a=2$ divides $2^{\nu(r)}\frac{s_i(r^f-1)}{q_i}$, by (\ref{nuMin}). Therefore,
    \begin{equation}\label{gamma}
    \gamma_{i}^{2^{\nu(r)}\frac{s_i(r^f-1)}{q_i}} =
    \left(\gamma_{i}^{c+1}\right)^{Y\frac{2^{\nu(r)}s_i(r^f-1)}{2q_i}} =
    \left(\beta_{i\sigma}^{s_i}\beta_{i\rho}^{-s_i}\right)^{2^{\nu(r)}Z}.
    \end{equation}
Let $j''\equiv j' \mod 2$ with $j''\in \{0,1\}$. Then $\Upsilon(\rho^{j''})=\Upsilon(\rho^{j'})$ and
$N_{\rho}^{j'}(w)=w^{j''}$. Therefore,
    \begin{equation}\label{Comm}
    \matriz{{rcl}
    [u_{\theta},u_{\rho}^{j'} u_{\sigma}^{j}] &=&
    [u_{\theta},u_{\rho}^{{j'}}]u_{\rho}^{{j'}}[u_{\theta},u_{\sigma}^{j}] u_{\rho}^{-{j'}}=
    \prod_{i=0}^n (\beta_{i\rho}^{-s_i})^{{j''}} (\beta_{i\sigma}^{-s_i})^{V(j)(-1)^{{j''}}} \\
    &=&
    \prod_{i=0}^n (\beta_{i\rho}^{-s_i})^{{j''}} (\beta_{i\sigma}^{-s_i})^{X\frac{c-1}{2^d}V(j)(-1)^{{j''}}} =
    \prod_{i=0}^n (\beta_{i\rho}^{-s_i})^{{j''}} (\beta_{i\sigma}^{-s_i})^{X\frac{c^j-1}{2^d}(-1)^{{j''}}}
    .}
    \end{equation}
Using (\ref{epstheta}), (\ref{gamma}) and (\ref{Comm}) we obtain
    \begin{equation}\label{FactorNC}
    \xi_0^{2^{\nu(r)}} = [u_{\theta},u_{\rho}^{j'}  u_{\sigma}^{j}]^{2^{\nu(r)}}
    \gamma_{\theta}^{2^{\nu(r)}\frac{r^f-1}{e}}
    = \left(\prod_{i=0}^n \beta_{i\rho}^{-s_i}\right)^{2^{\nu(r)}(Z+j'')}
      \left(\prod_{i=0}^n \beta_{i\sigma}^{s_i}\right)^{2^{\nu(r)}\left(Z-X\frac{c^j-1}{2^d}(-1)^{{j''}}\right)}.
    \end{equation}

We claim that $Z + j'' \equiv 0 \mod 2^{d-1}$. On the one hand $Y\equiv 1 \mod 2^{d-1}$. On the other hand, $\phi_q =
\rho^{j'} \sigma^j \eta_q$, with $\eta_q \in C_q$ and so $r^f\equiv (-1)^{j'}c^j\mod 2^{1+s+d}$. Hence $r^f\equiv
(-1)^{j'}=(-1)^{j''}\mod 2^d$ and therefore $Z+j''=Y\frac{r^f-1}{2}+j''\equiv \frac{(-1)^{j''}-1}{2}+j''\mod 2^{d-1}$.
Considering the two possible values of $j''\in \{0,1\}$ we have $\frac{(-1)^{j''}-1}{2}+j''=0$ and the claim follows.

From $d=v_2(c-1)$ one has $c\equiv 1 + 2^{d-1} \mod 2^d$ and hence $Y\equiv 1+2^{d-1} \mod 2^d$  and $r^f\equiv
(-1)^{j'}c^j \equiv (-1)^{j'}(1+j2^d)\mod 2^{1+s+d}$. Then
    $$\matriz{{rcl}
    \frac{Z+j''}{2^{d-1}} &=& \frac{Y(r^f-1)+2j''}{2^d} \equiv
    \frac{Y((-1)^{j''}(1+j2^d)-1)+2j''}{2^d} = \frac{Y(\frac{(-1)^{j''}-1}{2}+(-1)^{j''}j2^{d-1})+j''}{2^{d-1}} \\
    &\equiv&
    \frac{(1+2^{d-1})(-j''+(-1)^{j''}j2^{d-1})+j''}{2^{d-1}} =
    \frac{-j''-j''2^{d-1}+(-1)^{j''}j2^{d-1}+(-1)^{j''}j2^{2(d-1)}+j''}{2^{d-1}} \\
    &\equiv& -j''+(-1)^{j''}j \equiv j+j'' \equiv j+j' \mod 2.
    }$$
Using this, the equality $\beta_{i\rho}^{\frac{1-c}{2}} = \mu_i \beta_{i\sigma}$ and the fact that $\mu_i=\pm 1$ we
obtain
    $$\beta_{i\rho}^{-(Z+j'')}=\beta_{i\rho}^{-X\frac{c-1}{2^d}(Z+j'')}=
    \beta_{i\rho}^{-X\frac{c-1}{2}\frac{Z+j''}{2^{d-1}}} =
    \mu_i^{X\frac{Z+j''}{2^{d-1}}} \beta_{i\sigma}^{X\frac{Z+j''}{2^{d-1}}} =
    \mu_i^{j+j'} \beta_{i\sigma}^{X\frac{Z+j''}{2^{d-1}}}.$$
Combining this with (\ref{FactorNC}) we have
    $$\matriz{{rcl}
    \xi_0^{2^{\nu(r)}} &=&
 \prod_{i=0}^n \mu_i^{2^{\nu(r)}(j+j') s_i}
 \prod_{i=0}^n (\beta_{i\sigma}^{s_i})^{2^{\nu(r)}\left[Z-X\frac{c^j-1}{2^d}(-1)^{{j''}}+\frac{X(Z+j'')}{2^{d-1}}\right]} \\
  &=& \prod_{i=0}^n \mu_i^{2^{\nu(r)}(j+j') s_i}
  \prod_{i=0}^n (\beta_{i\sigma}^{s_i})^{2^{\nu(r)}\left[\frac{2^d Z+X(c^j-1)(-1)^{j''} +2X(Z+j'')}{2^d}\right]}.}$$
To finish the proof it is enough to show that the exponent of each $\beta_{i\sigma}$ in the previous expression is a
multiple of $2^{1+s}$.
Indeed, $2^d\equiv X(c-1) \mod 2^{1+s+d}$ and so
$$\matriz{{l}
    2^d Z+X(c^j-1)(-1)^{j''} +2X(Z+j'') \equiv ZX(c-1)-X(c^j-1)(-1)^{j''} +2X(Z+j'') = \\
    X(Y\frac{r^f-1}{2}(c+1)+(c^j-1)(-1)^{j''} +2j'') = X((r^f-1)Y\frac{c+1}{2}-c^j(-1)^{j''}+(-1)^{j''} +2j'') \equiv \\
    X({r^f-1}-c^j(-1)^{j''}+1) \equiv 0 \mod 2^{1+s+d}}$$
as required. This finishes the proof of the lemma in Case 2.
\end{proof}

We need the following Proposition from \cite{Jan}.

\begin{proposition}\label{Max}
For every odd prime $q\ne r$ not dividing $m$ let $d(q)=\min\{a,v_p(q-1)\}$. Then
\begin{enumerate}
\item $|c_0^{k_q} C/C^{p^{d(q)}}| \le |\theta_q^f C/C^{p^a}|$, and
\item the equality holds if $q\equiv 1 \mod p^a$ and $r$ is not
congruent with a $p$-th power modulo $q$.  There are infinitely
many primes $q$ satisfying these conditions.
\end{enumerate}
\end{proposition}

\begin{proof}
See Proposition 4.1 and Lemma 4.2 of \cite{Jan}.
\end{proof}

We are ready to prove the main result of the paper.

\begin{theorem}\label{localindex}
Let $K$ be an abelian number field, $p$ a prime and $r$ an odd prime. If either $\zeta_p\not\in K$ or $r\not\equiv
1 \mod p$ then $\beta_p(r)=0$. Assume otherwise that $\zeta_p\in K$ and $r\equiv 1 \mod p$, and use Notation~\ref{Not}
including the decomposition $\phi=\eta \rho^{j'}\sigma^j$ with $\eta\in B$.

\begin{enumerate}
\item Assume that $r$ does not divide $m$.
\begin{enumerate}
\item If $G/C$ is non-cyclic and $j\not\equiv j' \mod 2$ then $\beta_p(r)=1$.
\item Otherwise $\beta_p(r)=\max\{\nu(r),v_p(|\eta B^{p^{d(r)}}|)\}$, where $d(r)=\min\{a,v_p(r-1)\}$.
\end{enumerate}

\item Assume that $r$ divides $m$ and let $q_0$ be an odd prime not dividing $m$ such that $q_0\equiv 1 \mod p^a$ and
$r$ is not a $p$-th power modulo $q_0$. Let $\theta=\theta_{q_0}$ be a generator of the inertia group of $G_{q_0}$ at $r$.
\begin{enumerate}
\item If $G/C$ is non-cyclic, $j\not\equiv j' \mod 2$ and $\theta$ is not a square in $D$ then $\beta_p(r)=1$.

\item  Otherwise $\beta_p(r)=\max\{\nu(r), h, v_p(|\theta^f C^{p^a}|)\}$, where
$h=\max_{\Psi}\{v_p(|\Psi(\theta,\eta)|)\}$ as $\Psi$ runs over all skew pairings of $B$ over $\GEN{\zeta_{p^a}}$.
\end{enumerate}

\end{enumerate}
\end{theorem}

\begin{proof}
For simplicity we write $\beta(r)=\beta_p(r)$. We already explained why if either $\zeta_p\not\in K$ or $r\not\equiv 1
\mod p$ then $\beta_p(r)=0$. So in the remainder of the proof we assume that $\zeta_p\in K$ and $r\equiv 1 \mod p$, and
so $K$, $p$, and $r$ satisfy the condition mentioned at the beginning of the section.  It was also pointed out
earlier in this section that $p^{\beta(r)}$ is the $r$-local index of a crossed product algebra $A$ of the form
$A = (F(\zeta_q)/K, \alpha)$ with $q$
and $\alpha$ taking values in $\GEN{\zeta_{p^s}}$ or in $\GEN{\zeta_4}$. Moreover, since $p^{\nu(r)}$ is the $r$-local
index of the cyclic Schur algebra $(K(\zeta_r)/K, c_0, \zeta_{p^a})$ \cite{Jan}, we always have $\nu(r) \le \beta(r)$.

In case 1 one may assume that $q=r$, because $(F(\zeta_q)/K,\alpha)$ has $r$-local index $1$ for every $q\ne r$. Since
$\Gal(F(\zeta_r)/F)$ is the inertia group at $r$ in $G_r$, in this case one may assume that $\theta=\theta_r=c_0$. On
the contrary, in case 2, $q\ne r$, and $\theta=c_1^{s_1}\dots c_n^{s_n} \rho^{2s_{n+1}}$,
for some
$s_1,\dots,s_{n+1}$.

In cases (1.a) and (2.a), $G/C$ is non-cyclic and hence $p^a=2$. Then $\beta(r)\le 1$, by the Benard-Schacher Theorem,
and hence if $\nu(r)=1$ then $\beta(r)=1$. So assume that $\nu(r)=0$.  Furthermore, in case (2.a), $s_i$ is odd for
some $i\le n$, because $\theta \not\in D^2$.  Now we can use Corollary 5 to
produce a cyclotomic algebra $A'=(F(\zeta_q)/K,\alpha')$ so that $\xi_{\alpha}=-\xi_{\alpha'}$. Indeed, there is such
an algebra such that all the data associated to $\alpha$ are equal to the data for $A$, except for $\beta_{0\sigma}$,
in case (1.a), and $\beta_{k\sigma}$, case (2.a).  Using Lemma~\ref{xi0}
and the assumptions $\nu(r)=0$ and $j\not\equiv j'
\mod 2$, one has $\xi_{0,\alpha}=-\xi_{0,\alpha'}$ and $\Psi_{\alpha}=\Psi_{\alpha'}$. Thus
$\xi_{\alpha}=-\xi_{\alpha'}$, as claimed. This shows that $\beta(r)=1$ in cases (1.a) and (2.a).

In case (1.b),
$\xi=\xi_0 \Psi(c_0,\eta)$. By
Lemma~\ref{xi0}, $\xi_0$ has order dividing $p^{\nu(r)}$ in this case and, by Lemma~\ref{SPEx},
$\max\{|\Psi(\theta,\eta)|:\Psi \in S\}=|\eta B^{p^{d(r)}}|$, where $S$ is the set of skew pairings of $B_r$ with values in
$\GEN{p^a}$. Using this and $\nu(r)\le \beta(r)$ one deduces that $\beta(r)=\max\{\nu(r),v_p(|\eta B^{p^{d(r)}}|)\}$.

The formula for case (2.b) is obtained in a similar way using the equality $\xi = \xi_0 \Psi(\theta,\eta) \Psi(\theta,c_0^{s_0})$ and Lemmas~\ref{Maximums} and ~\ref{SPEx}.
\end{proof}

\end{document}